\documentclass[12pt]{article}

\usepackage{amsmath,amsfonts,amsthm,epsfig,latexsym,graphicx,amssymb}
\usepackage{amssymb,amsthm}

\newtheorem{Que}{Question}
\newtheorem{Cnj}{Conjecture}
\newtheorem{Thm}{Theorem}
\newtheorem{Corl}{Corollary}

\newtheorem{Lm}{Lemma}[section]
\newtheorem{Def}[Lm]{Definition}
\newtheorem{Prop}[Lm]{Proposition}
\newtheorem{Rem}[Lm]{Remark}
\newtheorem{Cor}[Lm]{Corollary}
\newtheorem{Cla}[Lm]{Claim}

\def\bdef{\begin{Def}}
\def\endef{\end{Def}}

\def\bthm{\begin{Thm}}
\def\bcnj{\begin{Cnj}}
\def\ethm{\end{Thm}}
\def\ecnj{\end{Cnj}}
\def\bque{\begin{Que}}
\def\eque{\end{Que}}

\def\brema{\begin{Rem}}
\def\erema{\end{Rem}}

\def\bpro{\begin{Prop}}
\def\epro{\end{Prop}}

\def\blm{\begin{Lm}}
\def\elm{\end{Lm}}

\def\bcorl{\begin{Corl}}
\def\ecorl{\end{Corl}}

\def\bcor{\begin{Cor}}
\def\ecor{\end{Cor}}

\def\brm{\begin{Rem}}
\def\erm{\end{Rem}}

\def\bcl{\begin{Cla}}
\def\ecl{\end{Cla}}

\def\bfig{\begin{picture}}
\def\efig{\end{picture}}

\def\beq{\begin{eqnarray}}
\def\eneq{\end{eqnarray}}

\def\beal{\begin{aligned}}
\def\enal{\end{aligned}}

 \def\NN{{\mathbb N}} 

 \def\RR{{\mathbb R}} \def\SS{{\mathbb S}}
\def\TT{{\mathbb T}}

\def\La{\Lambda}

\def\Om{\Omega}

\def\la{\lambda}

\def\ga{\gamma}
\def\const{\Theta}

\def\diffM{\mbox{{\rm Diff\,}}^1(M)}

  \def\cG{{\cal G}}  
  \def\cH{{\cal H}} \def\cN{{\cal N}} \def\cT{{\cal
T}}
\def\cC{{\cal C}}   \def\cO{{\cal O}} \def\cU{{\cal
U}}
    
    \def\cW{{\cal
W}}
   \def\cR{{\cal R}} 
  \def\cZ{{\cal Z}}

\def\st{\operatorname{s}}
\def\sst{\operatorname{ss}}

\def\ut{\operatorname{u}}
\def\uut{\operatorname{uu}}

\def\loc{\operatorname{loc}}

\title{Non-hyperbolic ergodic measures with large support}
\author{Ch. Bonatti, L. J. D\'\i az, and A. Gorodetski}

\makeatletter

\def\keywords#1{{\def\@thefnmark{\relax}\@footnotetext{#1}}}

\let\subjclass\keywords

\makeatother

\begin{document}

\maketitle

\keywords{Keywords: dominated splitting, homoclinic class, Lyapunov exponent, partial
hyperbolicity, heterodimenional cycle, support of an
invariant measure.}





\begin{abstract}
We prove that there is a residual subset $\mathcal{S}$ in
$\text{Diff\,}^1(M)$ such that, for every $f\in \mathcal{S}$, any
homoclinic class of $f$ with invariant one dimensional
central bundle containing saddles of different indices (i.e. with different dimensions of the stable invariant manifold)
coincides with the support of some invariant ergodic non-hyperbolic (one of the
Lyapunov exponents is equal to zero) measure of $f$.
\end{abstract}



\section{Introduction}

How to characterize the absence of uniform hyperbolicity? What
dynamical structures can not exist in the uniformly hyperbolic
setting but must be present in the complement?  In this paper, we
study how the non-hyperbolicity is detected in the ergodic level. Namely, we consider non-hyperbolic invariant measures as indication of non-hyperbolicity, and construct such measures with full support for partially hyperbolic homoclinic classes.

Some other candidates for the role of ``non-hyperbolic
structure" are the cycles (homoclinic tangencies and
heterodimensional cycles), super-exponential growth of the number
of periodic points \cite{K,BDF},  absence of shadowing property
\cite{BDT,YY,AD,Sa}, and non-existence of symbolic extensions
\cite{DN,A,DF}.

Recall that if $\mu$ is an ergodic measure of a diffeomorphism $f:
M\to M,\, \text{\rm dim}\, M={\mathbf{m}},$ then there is a set
$\Lambda$ of full $\mu$-measure and real numbers  $\chi^1_{\mu}\le
\chi^2_{\mu}\le \dots \le \chi^{\mathbf{m}}_{\mu}$ such that, for
every $x\in \Lambda$ and every non-zero vector $v\in T_xM$, one
has $\lim_{n\to \infty} 1/n \, \log || Df^n(v)|| = \chi^i_\mu$ for
some $i=1, \dots, \mathbf{m}$, see \cite{O, M2}. The number
$\chi^i_\mu$ is  the \emph{$i$-th Lyapunov exponent} of the
measure $\mu$. \bdef An ergodic invariant measure of a
diffeomorphism is called {\it non-hyperbolic} if at least one of
its Lyapunov exponents is equal to zero.
\endef

Inspired by Palis' density conjecture (\cite{Pa}), we stated the following (\cite{DG}):

\bcnj   In $\text{{\rm Diff\,}}^r(M), r\ge 1,$ there exists an
open and dense subset $\mathcal{U}\subset \text{{\rm
Diff\,}}^r(M)$ such that every diffeomorphism $f\in \mathcal{U}$
is either uniformly hyperbolic or has an ergodic non-hyperbolic
invariant measure. \ecnj

We consider the question of the existence of non-hyperbolic measures
 for homoclinic classes. The
{\emph{homoclinic class} of a hyperbolic periodic point $P$
of a diffeomorphism $f$, denoted by $H(P,f)$, is the closure of
the transverse intersections of the invariant manifolds of the
orbit of $P$, see \cite{N2}. A homoclinic class is a
{\emph{transitive set}} (there exists a dense orbit) where
periodic points also form a dense subset. Note that a homoclinic
class may fail to be hyperbolic and  may contain saddles of
different {\emph{$\st$-indices}} (dimension of the stable bundle).
This is precisely the setting we consider in this paper. In many important cases homoclinic classes are used to structure the global dynamics, playing a role
similar to that of the basic sets in the hyperbolic theory (in
fact, basic sets are a special case of homoclinic classes). In
\cite{DG} we show the following:

\medskip

\noindent {\bf Theorem.} {\em  There is a residual subset
$\mathcal{S}$  of $\text{\rm{Diff\,}}^1(M)$  such that, for every
$f\in \mathcal{S}$, any homoclinic class of $f$ containing saddles
of different $\st$-indices (dimension of the stable bundle)
contains also an uncountable support of an invariant ergodic
non-hyperbolic measure of $f$.\/}

\medskip

In view of this result, it is  natural to consider the following
question about the support of the ergodic non-hyperbolic measures
in non-hyperbolic homoclinic class.

\begin{Que}\label{question}
When does a non-hyperbolic homoclinic class  equal to the
support of a non-hyperbolic ergodic measure?
\end{Que}

The main result of this paper is the following (for an accurate formal statement see
Theorem~\ref{t.bdg}):

\medskip

\noindent {\bf Main result.} {\em For a $C^1$-generic diffeomorphism $f$, every homoclinic class $H$ with a 1-dimensional central direction and saddles of different indices is the support of an ergodic non-hyperbolic invariant measure $\mu$ of $f$,  $H=\text{\rm supp}\, \mu$.\/}

\medskip

In an opposite direction,  we would like to mention the results in
\cite{ABC}: for $C^1$-generic diffeomorphisms  generic measures
supported on  isolated homoclinic class are ergodic and hyperbolic
(all Lyapunov exponents are non-zero).

\medskip

The constructions in \cite{DG} translate to the setting of
homoclinic classes the arguments in \cite{GIKN} and \cite{KN} about the
existence of non-hyperbolic ergodic measures for  skew products
defined over the circle $\SS^1$ and the corresponding smooth realizations. There the circle
$\SS^1$ corresponds to the center direction and the fact that the
central direction is one-dimensional is essential. This condition
can be precisely formulated in terms of dominated splittings.

\begin{Def}[Dominated splitting and partial hyperbolicity]
Consider a diffeomorphism $f$ and a compact $f$-invariant set
$\La$. A $Df$-invariant splitting $T_{\La}M=E\oplus F$ over $\La$
is {\em{dominated\/}} if the fibers $E_x$ and $F_x$ of $E$ and $F$
have constant dimension and there is a constant $k\in \NN$ such
that
$$
\frac{||D_x f^k(u)||}{||D_xf^{k} (w) ||} <  \frac{1}{2},
$$
 for
every $x\in \La$ and every pair of unitary vectors $u\in E_x$ and
$w\in F_x$.

In some cases, we consider splittings with three bundles. A
$Df$-invariant splitting
$$
T_\La M=E\oplus F\oplus G
$$
over $\La$ is dominated if both
splittings $(E\oplus F)\oplus G$ and $E\oplus (F\oplus G)$
are dominated.

The dominated splitting $T_\La M=E\oplus F\oplus G$
is {\em{partially hyperbolic\/}} if
 $E$ and $G$ are both uniformly hyperbolic and at least one of them is not empty.
We say that $F$ is the {\em{central direction}} of the splitting.
\end{Def}

In \cite{DG} general  (non-hyperbolic) homoclinic classes are considered (i.e., there is
no hyperbolic-like assumptions).
A key step in \cite{DG} is to identify a partially
hyperbolic region $\La$ of the homoclinic class where the
(non-hyperbolic) central direction has dimension one: there is a
partially hyperbolic splitting $T_{\La}M=E\oplus F \oplus G$,
where $F$ is one-dimensional and non-hyperbolic and $E$ and $G$
are uniformly hyperbolic.

 In general, this partially hyperbolic region $\La$ is
properly contained in the whole homoclinic class, thus 
the support of the obtained measure is not the whole homoclinic class.
Prototypical examples of this case are the Derived from
Anosov diffeomorphisms defined on $\TT^3$ obtained via Hopf
bifurcations, see \cite{Ca,BV}, where the homoclinic classes are
the whole torus and its dominated splitting is of the form $T
\TT^3= E\oplus E^{\ut}$, where $E$ is two-dimensional, non-hyperbolic, and
undecomposable (it does not contain one-dimensional invariant
directions).
Thus in these cases the
arguments in \cite{DG} do not provide a non-hyperbolic measure
supported on the whole homoclinic class.

On the other hand, when the non-hyperbolic central bundle has a
one-dimensional invariant direction one can prove the following:

\bthm \label{t.bdg} Let $M$ be a closed manifold, $\dim M \ge 3$.
There is a residual subset $\cR$ of $\diffM$ such that for every
$f\in \cR$ and every homoclinic class $H(f)$ of $f$ such that
\begin{itemize}
\item
$H(f)$ has a dominated splitting $T_{H(f)}M=E\oplus F \oplus G$,
where $F$ has dimension one,
\item
$H(f)$ contains saddles of $\st$-indices $\dim (E)$ and
$\dim(E\oplus F)=\dim (E)+1$,
\end{itemize}
there is a non-hyperbolic ergodic $f$-invariant measure
$\mu_{H(f)}$ such that
$$
\mbox{\rm{supp}}\, \mu_{H(f)}= H(f).
$$
 \ethm

The simplest setting for our result is the one of homoclinic
classes having a partially hyperbolic splitting with a
one-dimensional central direction. There are two main sorts of
such homoclinic classes. On the one hand, there are the
non-hyperbolic homoclinic classes generated by unfolding
heterodimensional cycles, see \cite{D,DR}. On the other hand,
there are robustly transitive and non-hyperbolic
diffeomorphisms (in this case, the homoclinic class is the whole
ambient manifold). Recall that a diffeomorphism $f$ is
{\emph{robustly transitive}} if there is a neighborhood $\cU_f$ of
$f$ in $\diffM$ consisting of {\emph{transitive diffeomorphisms}}: every
$g\in \cU_f$ has a dense orbit in $M$.

Among the robustly transitive diffeomorphisms with one-dimensional
central direction, we mention the Derived from Anosov
diffeomorphisms (via saddle-node or fork bifurcations) in
\cite{M2}, the time-one maps of transitive Anosov vector fields,
see \cite{BD1}, and the perturbations of products of Anosov
diffeomorphisms and maps defined on the circle (skew products),
\cite{BD1,SW,RW}. For these cases, Theorem~\ref{t.bdg} just can be
read as follows:

\bcorl \label{c.dg} Let $M$ be a closed manifold, $\dim M \ge 3$.
There is a residual subset $\cR$ of $\diffM$ such that for every
$f\in \cR$ and every partially hyperbolic homoclinic
class
of $f$ having saddles of different indices and an one-dimensional
central bundle
 there is a non-hyperbolic ergodic $f$-invariant measure
whose support is the whole homoclinic class.
 \ecorl

We next discuss the interplay between the results above, the sorts
of dominated splittings of  homoclinic classes, and the occurrence
of homoclinic tangencies. Consider the subset $\mbox{HT}^1(M)\subset
\diffM$ of diffeomorphisms $f$ having a {\emph{homoclinic
tangency}} associated to some saddle (i.e., the invariant
manifolds of the saddle have some non-transverse intersection).
Define the (open) set of {\emph{diffeomorphisms far from
homoclinic tangencies} by
$$
\mbox{FT}^1(M) =\diffM\setminus \overline{\mbox{HT}^1(M)}.
$$

A recent result by Yang states the following dichotomy for
$C^1$-generic diffeomorphisms $f$ far from homoclinic tangencies
(i.e., diffeomorphisms $f$ in a residual subset of
$\mbox{FT}^1(M)$): any homoclinic class of $f$ is either
hyperbolic or it supports a non-hyperbolic ergodic measure, see
\cite{Y}. Note that, in
principle, there is the possibility for a generic diffeomorphism
$f$ to have  non-hyperbolic homoclinic classes such that all saddles have
the same $\st$-index. For homoclinic classes of diffeomorphisms in
$\mbox{FT}^1(M)$ containing saddles of different  indices, Theorem~\ref{t.bdg}
implies that there is  an ergodic
non-hyperbolic measure whose support is the whole homoclinic
class.

Finally, we state a refinement of Theorem~\ref{t.bdg} about the zero Lyapunov exponents of the
non-hyperbolic measure.}}
 Consider a homoclinic class $H(f)$ of $f$, its
\emph{$\st$-index (variation) interval} is the interval
$\st$-$\mbox{ind} (H(f))=[i,j]$, where $i$ and $j$ are the minimum
and the maximum of the \emph{$\st$-indices} of the periodic points
in $H(f)$. The homoclinic class $H(f)$ has \emph{index variation}
if $i<j$.

\bcorl There is a residual subset $\cR$ of $\diffM$ such that for
every diffeomorphism $f\in FT^1(M)\cap \cR$, every homoclinic
class $H(f)$ of $f$ with $\st$-index variation interval $[i,j]$,
and every $k\in [i,j)$ there is an ergodic measure $\mu^k_{H(f)}$
whose support is the whole homoclinic class $H(f)$ such that its
$k$-th Lyapunov exponent $\chi^{k}(\mu^k_{H(f)})$ is zero.
\label{c.gourmelon}
\ecorl

Indeed, fix $k\in [i,j)$. Theorem 1.1 in \cite{Go} claims that if $H(f)$ is a
homoclinic class such that $k,k+1$ belong to the index interval
$\st$-$\mbox{ind} (H(f))$  and there is no  dominated splitting
defined over $H(f)$ of the form $E\oplus F\oplus G$ with $\dim
(E)=k$ and $\dim(F)=1$, then there is a diffeomorphism $g$
arbitrarily $C^1$-close to $f$ having a homoclinic tangency
associated to a saddle of $\st$-index $k$.
Therefore, since $f \in \mbox{FT}^1(M)$,
the homoclinic class $H(f)$ has a
dominated splitting of the form $E\oplus F\oplus G$, where $\dim (E)=k$ and $F$ is one-dimensional.
On the other hand, by
 \cite{ABCDW}, as $f$ is a $C^1$-generic
diffeomorphism the
homoclinic class $H(f)$ contains saddles of
$\st$-indices $k,k+1\in [i,j]$.
Now we can  apply
Theorem~\ref{t.bdg} to get Corollary~\ref{c.gourmelon}.

\medskip

We close this introduction discussing the special case of robustly
non-hy\-per\-bo\-lic and transitive diffeomorphisms.  We point out
that the proof of our results makes a systematic use of  the
machinery recently developed about the dynamics and generic
properties of homoclinic classes of $C^1$ diffeomorphisms, see
Section \ref{s.generic} for details. Therefore our results are formulated for
$C^1$-generic diffeomorphisms.

On the other hand, the setting of robustly transitive and
non-hyperbolic diffeomorphisms is a natural framework for
considering Question~\ref{question} for open sets of
non-hyperbolic diffeomorphisms. Moreover, in a recent paper,
Nalsky gave examples of open sets of diffeomorphisms exhibiting
non-hyperbolic ergodic measures supported on the whole homoclinic
classes, \cite{Na}. In his case, one essentially has a skew
product over a hyperbolic diffeomorphism in a base and a circle as
a fiber. This example improves the construction
in \cite{GIKN} mentioned before.   Motivated by these results we want to formulate the partial case of the Conjecture
\ref{c.dg} that looks more approachable:

\bcnj\label{robustly} Let $M$ be a closed manifold with
$\dim(M)\ge 3$. Consider an open subset  $\cT(M)$ of $\text{{\rm
Diff\,}}^1(M)$ consisting of non-hyperbolic transitive
diffemorphisms $f$ such that there is a $Df$-invariant dominated splitting  $TM=E\oplus F\oplus
G$, where $F$ has dimension one and is non-hyperbolic.
Then there is an open and dense subset $\cN(M)$ of $\cT(M)$ consisting of
diffeomorphisms $f$ having a non-hyperbolic ergodic measure
$\mu_f$
with full support, $\mbox{\rm{supp}}\, \mu_f= M$.
 \ecnj

 Notice that the approach used in this paper cannot be used straightforwardly to settle this conjecture. Indeed, our construction involves measures supported on periodic orbits, and therefore periodic orbits are dense in the support of the measure that we obtain. It is an open question whether periodic orbits are dense in a nonwandering set for an open and dense subset of the set of robustly transitive diffeomorphisms;  $C^1$-Closing Lemma \cite{Pu} implies this property only for a residual set of $C^1$-diffeomorphisms.

\medskip

This paper is organized as follows. In Section~\ref{s.ergodicmeasures} we develop the idea initially suggested
by Ilyashenko in \cite{GIKN} for constructing ergodic invariant measures as limit of measures supported on periodic orbits.
In Section~\ref{s.generic} we state properties of homoclinic classes of $C^1$-generic diffeomorphisms and recall results about
heterodimensional cycles. Finally, in Section~\ref{s.final} we construct a collection of atomic measures supported on
periodic orbits whose limit is a non-hyperbolic ergodic measure with full support, thus proving our main result.

\subsection*{Acknowledgments}

The authors thank the financial support and hospitality of ICTP.
We also express our gratitude to the organizers
of the conferences ``School and Workshop on Dynamical Systems", Trieste (2008), for creating a productive atmosphere for developing this paper.
Also we are grateful to the organizers of the  ``International Workshop on Global Dynamics Beyond Uniform Hyperbolicity'', Beijing (2009),
where this paper was completed.
L.J.D. thanks the financial support and warm hospitality of
UC Irvine during his stay for preparing this paper.
Ch.B. was partially supported by  ANR project ``DynNonHyp'' BLAN08-2-313375, L.J.D. by CNPq, Faperj (Cientista do Nosso Estado and Infraestrutura), and PRONEX, and A.G.  by grant NSF DMS-0901627.

\section{Non-hyperbolic ergodic measures as limits of atomic measures}
\label{s.ergodicmeasures}

The approach suggested in \cite{GIKN} allows to construct
non-hyperbolic ergodic invariant measures as limits of measures
supported on special sequences of periodic points. Here we show
that the support of the limit measure constructed in this way is
a ``topological limit" of the sequence of periodic points, see
Proposition~\ref{p.firstmainproposition} for details.

\subsection{Ergodicity, invariant direction fields, and Lyapunov exponents}
\label{ss.direction} Consider a diffeomorphism $f:M\to M$ and a closed invariant set $\Delta\subset M$. Assume that there is a $Df$-invariant
continuous direction field $E=(E_x)_{x\in \Delta}$ in $\Delta$. Then for
every invariant measure $\mu$ whose support is contained in $\Delta$
one of the Lyapunov exponents of $\mu$ is associated to $E$
(denote it by $\chi^E$). Namely, for $\mu$-a.e. $x\in M$ and for
every non-zero vector $v\in E_x\subset T_xM$,
$$
\lim_{n\to \infty}\frac{1}{n}\log|Df^n(v)|=\chi^E(\mu).
$$

In this section, by convergence of a sequence of measures we mean
$*$-weak convergence: $\mu_n$ converges to  $\mu$, if for any
continuous function $\varphi\colon M\to \mathbb{R}$ it holds
$$
\int \varphi \, d\mu_n \to \int \varphi \, d\mu, \qquad \mbox{as }
n\to\infty.
$$

To obtain non-hyperbolic measures we will use the following
statement.

\blm [\cite{DG}, Lemma 2.1, or \cite{GIKN}, Lemma 1] Let
diffeomorphism $f:M\to M$ have an invariant continuous direction
field $E$ in an $f$-invariant closed set $\Delta\subset M$.  Let
$\mu_n$ and $\mu$ be ergodic probability measures with supports in
$\Delta$, and $\mu_n\to \mu$ as $n\to \infty$. Then
$\chi^E(\mu_n)\to \chi^E(\mu)$. \label{l.limit} \elm

\subsection{Sufficient conditions for ergodicity and convergence}
\label{ss.sufficientergodicity}

In this section we state
an improved version of Lemma 2 from \cite{GIKN}. 
We need the following definition.

\bdef [$n$-measure] Let $\mathfrak{G}$ be a continuous map of a
metric compact space $\mathfrak{Q}$ into itself. A {\em
$n$-measure of the point $x_0$\/} is an atomic measure uniformly
distributed on $n$ subsequent iterations of the point $x_0$ under
the map $\mathfrak{G}$:
$$
\nu_n(x_0) = \frac{1}{n} \sum\limits_{i=0}^{n-1}
\delta_{\mathfrak{G}^i(x_0)},
$$
where $\delta_x$ is  $\delta$-measure supported at point $x$.
\endef

\blm\label{ergodicitylemma} Let
$\{X_n\}$ be a sequence of periodic orbits with increasing periods
$\pi(X_n)$ of a continuous map $\mathfrak{G}$ of a compact metric
space $\mathfrak{Q}$ into itself. For each $n$, let $\mu_n$ be the
probability atomic measure uniformly distributed on the orbit
$X_n$. 

Assume that for every $\varepsilon>0$ and every $m\in \mathbb{N}$
there exits a subset $\widetilde{X}_{m,\varepsilon} \subset X_m$
such that the following holds. For each continuous function
$\varphi$ on $\mathfrak{Q}$ there exists $N = N(\varepsilon,
\varphi) \in \Bbb{N}$ such that  for all $m>N$
   the following conditions are satisfied:
\begin{enumerate}
\item $ \mu_m(\widetilde{X}_{m,\varepsilon}) =\dfrac{\# \widetilde{X}_{m,\varepsilon}}{\# X_m} > 1 -
\varepsilon$ ($\# X$ the cardinality of the finite set $X$) and
\item for any  $n$, such that $m>n\ge N$, and for all $x \in
\widetilde{X}_{m,\varepsilon}$ it holds
$$
\left| \int\varphi \, d\nu_{\pi(X_n)}(x) - \int\varphi \, d\mu_n
\right| < \varepsilon.
$$
\end{enumerate}
Then the sequence of atomic measures $\{\mu_n\}$
has a limit, and  the limit
measure is ergodic. \elm

\begin{proof}
Lemma 2 from \cite{GIKN} claims that under conditions of Lemma \ref{ergodicitylemma} every limit point $\mu$  of the sequence $\{\mu_n\}$ is an
ergodic measure. Therefore we just need to show that the sequence of measures $\{\mu_n\}$ has a
$*$-weak limit. In order to do that it is enough to show that for any continuous
function $\varphi$ on $\mathfrak{Q}$ the sequence $\left\{\int
\varphi d\mu_n\right\}$ is Cauchy and thus converges. Indeed, since $\{\mu_n\}$ has a
convergent subsequence, $\mu_{n_k}\to \mu$ as $k\to \infty$, this
implies that $\int \varphi d\mu_n\to \int \varphi d\mu$ and
therefore $\mu_n\to \mu$ as $n\to \infty$.

Fix a continuous function $\varphi$ on $\mathfrak{Q}$. 
It is enough to show that the sequence $\left\{\int
\varphi d\mu_n\right\}$ is Cauchy. 
Since
$\mathfrak{Q}$ is compact, $\varphi$ is bounded, $|\varphi|<M$ for
some $M>0$. 
Fix small $\varepsilon>0$ and
consider the sets $\widetilde{X}_{m,\varepsilon} \subset X_m$ such
that for all  $m>N = N(\varepsilon, \varphi)$ properties (1) and
(2) hold. Now, for a given $m>n>N(\varepsilon, \varphi)$, we have
$$
\begin{array}{ll}
    \displaystyle{\sum_{x\in X_m}} \int\varphi d\nu_{\pi(X_n)}(x)&=
    \displaystyle{\sum_{x\in
    X_m}}\left(\dfrac{1}{\pi(X_n)}\, \displaystyle{\sum_{i=0}^{\pi(X_n)-1}}\varphi(\frak{G}^i(x))\right)=
    \\\\
     &=\dfrac{1}{\pi(X_n)}\,
     \displaystyle{\sum_{x\in
    X_m}\sum_{i=0}^{\pi(X_n)-1}}\varphi(\frak{G}^i(x))=\dfrac{1}{\pi(X_n)}\,\displaystyle{\sum_{x\in
    X_m}}\pi(X_n)\, \varphi(x)= \\\\
    &= \displaystyle{\sum_{x\in
    X_m}}\varphi(x)=\pi(X_m)\,\int\varphi d\mu_m.
\end{array}
$$
Therefore we can estimate
$$
\begin{array}{ll}
   \left|\int\varphi d\mu_m-\int\varphi
    d\mu_n\right|&= \dfrac{1}{\pi(X_m)}\left|\pi(X_m)\int\varphi d\mu_m-\pi(X_m)\int\varphi
    d\mu_n\right|=\\\\
    &=\dfrac{1}{\pi(X_m)}\left| \displaystyle{\sum_{x\in X_m}}\int\varphi d\nu_{\pi(X_n)}(x)-\pi(X_m)\int\varphi
    d\mu_n\right|\le\\\\
    &\le \dfrac{1}{\pi(X_m)}\left|\displaystyle{\sum_{x\in \widetilde{X}_{m, \varepsilon}}}\left(\int\varphi d\nu_{\pi(X_n)}(x)-\int\varphi
    d\mu_n\right)\right|+\\\\ &\quad+\dfrac{1}{\pi(X_m)}\left| \displaystyle{\sum_{x\not\in \widetilde{X}_{m, \varepsilon}}}
    \left(\int\varphi d\nu_{\pi(X_n)}(x)-\int\varphi
    d\mu_n\right)\right|.
\end{array}
$$
Finally, from condition (2) in the lemma, $|\varphi|<M$, and
condition (1), we have
$$
\left|\int\varphi d\mu_m-\int\varphi
    d\mu_n\right|
    \le
    \dfrac{1}{\pi(X_m)}\big(\pi(X_m)\varepsilon+2\,M\,\pi(X_m)\varepsilon\big)=(1+2\,M)\varepsilon.
$$
So the sequence $\left\{\int \varphi d\mu_n\right\}$ is Cauchy,
and Lemma~\ref{ergodicitylemma} is proved.
\end{proof}

\subsection{Sufficient conditions for existence of an invariant non-hyperbolic measure with large support}
\label{ss.sufficientnonhyperbolic}

In this section, we state an improved version of Proposition 2.5
from \cite{DG}, where we add a description of the support of the
resulting limit measure. We need the following definition.

\bdef [Good approximations] \label{d.goodaprox} A periodic orbit
$Y$ of a map $\mathfrak{G}$ of a compact metric
space $\mathfrak{Q}$ into itself is a $(\gamma, \varkappa)$-good approximation of
the periodic orbit $X$ of $\mathfrak{G}$ if the following holds.
\begin{itemize}
  \item There exists a subset $\Gamma$ of \,$Y$
and a projection $\rho \colon \Gamma \to X$ such that
$$
\text{\rm dist} (\mathfrak{G}^j(y), \mathfrak{G}^j(\rho(y)))<\gamma,
$$
for every $y\in \Gamma$ and every $j=0,1,\dots, \pi(X)-1$;
 \item $\dfrac{\#\Gamma}{\#Y}\ge \varkappa$; 
  \item $\#\rho^{-1}(x)$ is the same for all $x\in X$.
\end{itemize}
\endef

\blm\label{l.25} Let
$\{X_n\}$ be a sequence of periodic orbits with increasing periods
$\pi(X_n)$ of a continuous map $\mathfrak{G}$ of a compact metric
space $\mathfrak{Q}$ into itself. For each $n$, let $\mu_n$ be the
probability atomic measure uniformly distributed on the orbit
$X_n$.

Assume that there are sequences of numbers $\{\gamma_n\}_{n\in \mathbb{N}}$, $\gamma_n>0$, and $\{\varkappa_n\}_{n\in \mathbb{N}}$, $\varkappa_n\in (0, 1]$, such that
\begin{enumerate}
  \item for each $n\in \mathbb{N}$ the orbit $X_{n+1}$ is a $(\gamma_n, \varkappa_n)$-good approximation of $X_n$;
  \item $\sum_{n=1}^{\infty}\gamma_n<\infty$;
  \item $\prod_{n=1}^{\infty}\varkappa_n\in (0, 1]$. 
\end{enumerate}
Then the sequence of atomic measures $\{\mu_n\}$ has a limit $\mu$. The limit measure $\mu$ is ergodic, and
$$\text{\rm supp}\ \mu
=\bigcap_{k=1}^{\infty} \left(\,
\overline{\bigcup_{l=k}^{\infty}X_l }\, \right) \equiv \mathbf{X}.$$ In other words,
$\text{\rm supp}\ \mu$ is a topological limit of the sequence
of 
orbits $X_n$,
$$
\text{\rm supp}\ \mu =\{y\in \mathfrak{Q}\ |\ \exists m_i\to \infty \ \text{\rm and} \ x_i\in X_{m_i}\ \text{\rm such that} \ \lim_{i\to \infty}x_i=y\}.
$$
 \elm
\begin{proof}
 Let us  check first that the conditions of the Lemma \ref{ergodicitylemma} are
  satisfied by the sequence of orbits $X_n$.

Take arbitrary $\varepsilon > 0$ and continuous map $\varphi
\colon \mathfrak{Q} \to \mathbb{R}$. By assumption 1), for orbits $\{X_n\}$ a
sequence of subsets $\widetilde{X}_n\subset X_n$ and a sequence of
projections $\rho_n: {\widetilde X}_{n+1} \to X_n$ are defined
such that:

\begin{equation}
\label{e.varkappa} \prod\limits_{n=1}^{\infty} \dfrac{\#
{\widetilde X}_{n+1} } {\# X_{n+1}}\ge \prod\limits_{n=1}^{\infty}
\varkappa_n 
     > 0.
     \end{equation}

Choose $\delta = \delta (\varepsilon, \varphi)$ such that:
$$
\omega_{\delta} (\varphi) :=
    \sup\limits_{\text{\rm dist}(x,y) < \delta} |\varphi(x) - \varphi(y)| <
    \varepsilon.
$$

By assumption 2), we have  $\sum_{n=1}^{\infty}\gamma_n<\infty$.
Choose $N = N (\varepsilon, \varphi)$ such that the following
holds:
\begin{equation}
\sum\limits_N^{\infty} \gamma_k < \delta (\varepsilon, \varphi)
\qquad \mbox{and} \qquad \prod\limits_N^{\infty} \varkappa_k  > 1
- \varepsilon.
\label{e.star}
\end{equation}

Since the number of points in a pre-image for projections $\rho_n$
does not depend on a point in the image, a set ${\widetilde
X}_{m,\varepsilon} \subset X_m$
 where the total projection

$$
\rho_{m,N} = \rho_{m-1} \circ \dots \circ \rho_N\colon {\widetilde
X}_{m,\varepsilon}\to X_N
$$
is defined
 contains most of the orbit $X_m$:
$$
\frac { \# {\widetilde X}_{m,\varepsilon} } {\# X_m} \ge
    \prod \limits_{k=N}^{m-1} \varkappa_k \geq
    \prod\limits_N^{\infty} \varkappa_k >
    1 - \varepsilon.
$$
This implies 1) in Lemma~\ref{ergodicitylemma}.

Take arbitrary $m$ and $n$ with $m > n > N (\varepsilon,
\varphi)$. By construction, on the set ${\widetilde
X}_{m,\varepsilon}$ the total projection $\rho_{m,n} = \rho_{m-1}
\circ \dots \circ \rho_n$ is defined and
 for every point $x$ from the set
 $\widetilde{X}_{m,\varepsilon}\subset X_m$ from the first part of equation
\eqref{e.star} we have

\begin{equation}
\text{\rm dist} (\mathfrak{G}^j(x), \mathfrak{G}^j(\rho_{m,\,n}(x)))<\delta(\varepsilon,
\varphi), \qquad \text{for all $j=0, 1, \ldots, \pi(X_n)-1$.}
\label{e.2star}
\end{equation}
 Hence for $x \in {\widetilde X}_{m,\varepsilon}$ we have:
$$
\left| \int\varphi \, d\nu_{\pi(X_n)}(x) - \int\varphi \, d\mu_n
\right| <
    \omega_{\delta} (\varphi) <
    \varepsilon.
$$
Thus all conditions of Lemma~\ref{ergodicitylemma} are verified.
Therefore the sequence $\{\mu_n\}$ has a limit $\mu$, and the limit measure $\mu$ is ergodic.

 Now we need
the following lemma. Denote by $U_\delta(x)$ the ball of radius
$\delta$ centered at $x$.

\blm\label{l.m}  Set $r_n=\sum\limits_{k=n}^{\infty}
  \gamma_k $. For every point $x\in X_n$ 
one has $\mu(\overline{U_{r_n}(x)})>0$. \elm

\begin{proof}
Take any point $x \in X_n$. Notice
 that in its $\gamma_{n}$--neighborhood  there are at least
$\frac{\# {\widetilde X}_{n+1, \varepsilon}}{\pi(X_n)}$ points of
the orbit of $X_{n+1}$, where
\begin{equation}
\frac{\# {\widetilde X}_{n+1, \varepsilon}}{\pi(X_{n})} = \frac{\# {\widetilde X}_{n+1, \varepsilon}}{\pi(
X_{n+1})} \,\, \frac{\pi ( {X}_{n+1})}{ \pi( X_{n})} \ge \varkappa_{n}\,
 \frac{\pi(X_{n+1})}{\pi(X_n)}\equiv \bar \varkappa_{n}.
\label{e.3star}
\end{equation}
Therefore,
\begin{equation}
\label{e.medida}
 \mu_{n+1}(U_{\gamma_{n}}(x))\ge
\frac{\varkappa_{n}\,
 \frac{\pi(X_{n+1})}{\pi(X_n)}}{\pi(X_{n+1})}=
\varkappa_{n}\, \frac{1}{\pi(X_n)}= \varkappa_{n}\, \mu_n (\{x\}).
\end{equation}
In the neighborhood $U_{\gamma_n}(x)$ there are $p$
different points $x_1,\dots, x_p$ of the orbit of $X_{n+1}$, where
$p\ge \bar\varkappa_n$. Notice  that, by equation~\eqref{e.3star}, for each of these points $x_i$ the pre-image $\rho^{-1}_{n+2}(x_i)$ consists of at least
$$
\frac{\# \widetilde{X}_{n+2, \varepsilon}}{\pi(X_{n+1})}\ge \varkappa_{n+1}\frac{\pi(X_{n+2})}{\pi(X_{n+1})}
$$
points. The sets $\rho^{-1}_{n+2}(x_i)$ and $\rho^{-1}_{n+2}(x_j)$ are disjoint for $i\ne j$. Therefore we have
$$
\begin{array}{ll}
\mu_{n+2}(U_{\gamma_{n+1}+\gamma_n} (x)) &\ge \left[p\,
\varkappa_{n+1}\,
 \dfrac{\pi(X_{n+2})}{\pi(X_{n+1})}\right]\dfrac{1}{\pi(X_{n+2})} \ge \\\\ &\ge
 \left(
\varkappa_{n}\,
\dfrac{\pi(X_{n+1})}{\pi(X_n)}\right)\,\varkappa_{n+1}\,
 \dfrac{1}{\pi(X_{n+1})}=\displaystyle{\dfrac{\varkappa_n\,
 \varkappa_{n+1}}{\pi(X_{n})}}=\\\\&=
\varkappa_n\,
 \varkappa_{n+1}\, \mu_n(\{x\}).
 \end{array}
$$
Thus arguing inductively we have, for every $x\in X_n$,
$$
\mu_{n+\ell}(U_{\gamma_{n+\ell}+\cdots+
\gamma_{n+1}+\gamma_{n}}(x)) \ge (\varkappa_{n+\ell}\, \cdots\,
\varkappa_{n+1}\, \varkappa_{n}) \, \mu_n (\{x\}).
$$
Taking a limit and taking into account the condition {\it 3} in Lemma~\ref{l.25}, we have:
$$
\mu(\overline{U_{r_n} (x)}) \ge \left( \prod_{k=n}^{\infty}
\varkappa_{k}\right)\,
 \mu_n(\{x\}) >0, \quad \text{where $r_n=\sum\limits_{k=n}^{\infty}
\gamma_k$.} $$ Therefore, Lemma \ref{l.m} holds. \end{proof}

Now let us show that the support of $\mu$ is the topological limit
 $\mathbf{X}$ of the sequence of orbits $X_n$. It is a general fact that for any convergent
sequence of measures the support of the limit is a subset of the
topological limit of the sequence of supports, therefore
$\text{\rm supp}\ \mu\subseteq \mathbf{X}$. We need to show that
in our case we also have $\mathbf{X}\subseteq\text{\rm supp}\
\mu$.

Take any point $z\in \mathbf{X}$. It is enough to show that arbitrary small $\varepsilon>0$ one has $\mu(U_{\varepsilon}(z))>0$. Choose $n'\in \mathbb{N}$ large enough to guarantee that
$$
X_{n'}\cap U_{\varepsilon/2}(z)\ne \emptyset\ \ \ \text{\rm and} \ \ \ r_{n'}<\varepsilon/2.
$$
Take any $z'\in X_{n'}\cap U_{\varepsilon/2}(z)$. 
We have $U_{\varepsilon}(z)\supset \overline{U_{r_{n'}}(z')}$, and due to Lemma \ref{l.m} $\mu(\overline{U_{r_{n'}}(z')})>0$. Hence $\mu(U_{\varepsilon}(z))>0$, and therefore $z\in \text{\rm supp}\ \mu$.

This completes the proof of Lemma \ref{l.25}.
\end{proof}

\bpro \label{p.firstmainproposition} Assume that a diffeomorphism
$f:M\to M$ has the following properties:
\begin{description}
\item{{\bf 1)}}
there exists an $f$-invariant closed set $\Delta\subset M$ such that $f$ has
  an invariant continuous direction field $E$ in $\Delta$;
\item{{\bf 2)}}
there exists a sequence of periodic orbits
  $\{X_n\}_{n=1}^{\infty}, \, X_n\subset \Delta,$ of $f$
whose periods $\pi(X_n)$
  tend to infinity as $n\to \infty$.
\end{description}
   Denote by $\chi^E(X)$ the Lyapunov exponent along the periodic orbit $X$ with
respect to the invariant direction field $E$.
\begin{description}
\item{{\bf 3)}}
There exists a sequence of numbers $\{\gamma_{n}\}_{n=1}^{\infty},
\gamma_n>0$, and a constant $C>0$ such that for each $n$ the orbit
$X_{n+1}$ is a $(\gamma_n, 1-C\,|\chi^E(X_n)|)$-good approximation
of the orbit $X_{n}$;
\item{{\bf 4)}}
let $d_n$ be the minimal distance between the points of the orbit
$X_n$, then
$$
\displaystyle{\gamma_n< \frac{\min_{1\le i\le n} d_i}{3\cdot 2^n}};
$$
\item{{\bf 5)}} there exits a constant $\xi\in(0,1)$ such that for every $n$
$$|\chi^E(X_{n+1})|< \xi\, |\chi^E(X_n)|.$$
\end{description}
Then the sequence of atomic measures $\mu_n$
supported on periodic orbits $X_n$ has a ($*$-weak) limit $\mu$. The
measure $\mu$ is a non-hyperbolic ergodic invariant measure of $f$, and $\text{\rm supp}\ \mu$ is uncountable. Moreover, $\text{\rm supp}\ \mu$ is a topological limit of the sequence
of
orbits $X_n$, i.e.
$$\text{\rm supp}\ \mu =\bigcap_{k=1}^{\infty} \left(\,
\overline{\bigcup_{l=k}^{\infty}X_l }\, \right).$$
  \epro

\begin{proof}
Assumptions {\bf 1)}--{\bf 5)} of Proposition
\ref{p.firstmainproposition} imply that conditions of Lemma
\ref{l.25} are satisfied.  Therefore there
exists a limit $\lim_{n\to \infty}\mu_n=\mu$,  the measure
$\mu$ is an ergodic invariant measure of $f$, and $\text{\rm supp}\ \mu$ is a topological limit of the sequence
of orbits $X_n$.  Due to Lemma
\ref{l.limit} and property {\bf 5)}, measure $\mu$ is
non-hyperbolic. Finally, Proposition 2.5 from \cite{DG} implies that $\text{\rm supp}\ \mu$ is uncountable.
\end{proof}

\section{Homoclinic classes: $C^1$-generic properties and generation of cycles}\label{s.generic}


In this section we recall some known results on homoclinic classes and heterodimensional cycles that will be used later in our construction.


\subsection{Generic properties}
\label{ss.generic}

Here we state properties of homoclinic classes of
$C^1$-generic diffeomorphisms. There is a residual subset
${\cal{G}}$ of $\diffM$ such that every diffeomorphism $f\in \cG$
satisfies properties {\bf R1)}--{\bf R5)}  below.
\begin{description}
\item{\bf R1)}
Every homoclinic class of $f\in \cG$ depends continuously on $f$ in Hausdorff metric,
see \cite{CMP}. Moreover, if the homoclinic class
 contains saddles of indices
$a$ and $b$, $a<b$, it also contains saddles of index $c$, for
every $c\in (a,b)\cap\NN$. See \cite[Theorem 1]{ABCDW}.
\end{description}

Consider a hyperbolic periodic point $P_f$ of a diffeomorphism
$f$. It is well known that there are open neighborhoods $U$ of
$P_f$ in the  manifold and $\mathcal{U}$ of $f$ in $\diffM$ such
that every $g\in \cU$ has a unique hyperbolic periodic point $P_g$
of the same period as $P_f$ in $U$. The point $P_g$ is called the
{\it continuation} of $P_f$.

Recall that two saddles $P_f$ and $Q_f$ are {\emph{homoclinically related}} if their
invariant manifolds $W^{\st}(P_f)$ and  $W^{\ut}(Q_f)$ and
$W^{\ut}(P_f)$ and $W^{\st}(Q_f)$
have non-empty transverse intersections. In this case, the 
homoclinic classes of $P_f$ and of $Q_f$ coincide and the saddles have the same index.
Moreover, the  continuations $P_g$ and $Q_g$ are also homoclinically related
for all $g$ close to $f$.

\bdef[Persistently linked saddles] Consider a pair of hyperbolic saddles $P_f$ and $Q_f$ whose
continuations are defined for every $f$ in an open set\, $\cU$ of
$\diffM$. The saddles $P_f$
and $Q_f$ are {\em{persistently linked in $\cU$}} if there is a
residual subset $\cR$ of $\cU$ such that $H(P_f,f)=H(Q_f,f)$ for
all $f\in \cU\cap \cR$. \label{d.linked}
\endef

\begin{description}
\item{\bf R2)}
Given any $f\in \cG$ and any pair of saddles $P_f$ and $Q_f$ of
$f$, there is a neighborhood $\cU_f$ of $f$ such that either $P_f$
and $Q_f$  are persistently linked in $\cU_f$ or $H(P_g,g)\cap
H(Q_g,g)=\emptyset$ for all $g\in \cU_f\cap \cG$. 
Moreover, If the saddles $P_f$ and $Q_f$ have the same index
they are homoclinically related.
See \cite[Lemma
2.1]{ABCDW}.
\end{description}

\bdef [Saddle with real multipliers]
Let $P$ be a periodic point of period $\pi(P)$ of a 
 diffeomorphism $f$. We say that
$P$ has {\emph{real multipliers}} if every eigenvalue $\lambda$ of
$Df^{\pi(P)}(P)$ is real and has multiplicity one, and two
different eigenvalues of $Df^{\pi(P)}(P)$ have different absolute
values. We order the eigenvalues of $Df^{\pi(P)}(P)$ in increasing
ordering according their absolute values $|\la^1(P)|<\cdots
<|\la^{\mathbf{m}}(P)|$ and say that $\la^k(P)$ is the
\emph{$k$-th multiplier of $P$}. \label{d.real}
\endef
\bdef [$\st$- and $\ut$-index]
Let $P$ be a hyperbolic saddle. 
 By {\rm $\st$-index} of $P$ we mean the number of multipliers of $P$ with absolute value less than one, and  by {\rm $\ut$-index}  the number of multipliers of $P$ with absolute value greater than one.
\endef

Let $P$ be a saddle  with real multipliers, and suppose that  {\rm $\st$-index}$(P)=t+1$.
Consider the bundle $E^{\sst}\subset T_PM$ corresponding to the
first $t$ contracting eigenvalues of $P$ and the {\emph{strong
stable manifold}} $W^{\sst}(P)$ of $P$ (the only $f$-invariant
manifold of dimension $t$ tangent to the strong stable direction
$E^{\sst}$).

\bdef[$\st$- and $\ut$-biaccumulation]\label{biacc} A hyperbolic
saddle $P$  with real multipliers  is {\emph{$\st$-biaccumulated}}
(by its homoclinic points) if both connected components of
$W^{\st}_{\loc}(P)\setminus W^{\sst}_{\loc}(P)$ contain transverse
homoclinic points of $P$. Define also {\emph{$\ut$-biaccumulation
by homoclinic points}} in a similar way.
\endef

Note that $\st$- and $\ut$-biaccumulation are open
properties.

\medskip

Given a saddle $P$, we denote by $\mbox{Per}_{\RR}(H(P,f))$ the
saddles homoclinically related to $P$ having real multipliers.
 Clearly,
 $\mbox{Per}_{\RR}(H(P,f))\subset H(P,f)$.

\begin{description}
\item{\bf R3)}
For every diffeomorphism $f\in \cG$ and every saddle $P$ of $f$
whose homoclinic class is non-trivial the set
$\mbox{Per}_{\RR}(H(P,f))$ is dense in the whole homoclinic class
$H(P,f)$.
Moreover, the saddles from $\mbox{Per}_{\RR}(H(P,f))$
which are $\st$- (or $\ut$-) biaccumulated are also dense in $H(P,f)$.
See \cite[Proposition 2.3]{ABCDW} and \cite[Lemma 3.4]{DG}.

\item{\bf R4)}
Consider $f\in \cG$ and a saddle $P$ of $f$ whose homoclinic class
is non-trivial. For every $\epsilon>0$ and every $k=1,\dots,
{\bf{m}}$, one has that
$$
\overline{\Big(
 \{S\in \mbox{Per}_{\RR}(H(P,f))\colon
 |\chi^k(S)-\chi^k(P)|<\epsilon\}\Big)}=H(P,f).
 $$
Assume now that $\chi^k(P)<0$ and suppose that $H(P,f)$ contains a
saddle $Q$ with $\chi^k(Q)>0$. Then one has that
$$
\overline{ \Big(
 \{S\in \mbox{Per}_{\RR}(H(P,f))\colon \chi^k(S)\in
 (-\epsilon,0)\}\Big)}=H(P,f).
 $$
 Moreover, we can additionally assume that the dense subsets of
$H(P,f)$ above consist of saddles with the biaccumulation
properties. For details see
 \cite{ABCDW,BDF,DG}.

\item{\bf R5)} We will actually need a slightly stronger property.

\blm\label{l.perdense} Consider any $f\in \cG$ and any hyperbolic saddle $P_f$
whose homoclinic class is non-trivial. For every $\epsilon>0$ and $k=1,\dots,
{\bf{m}}$, the homoclinic class $H(P_f,f)$ of $f$ contains a
saddle $Y_f\in H(P_f,f)$ with real multipliers
such that $Y_f$ is homoclinically related to
$P_f$, its orbit is $\epsilon$-dense in the homoclinic class
$H(P_f,f)$, and the Lyapunov exponent
$\chi^k(Y_f)$ is $\epsilon$-close to $\chi^k(P_f)$. \elm

\begin{proof}
Indeed, for a given $\varepsilon>0$, by {\bf R4)} there exists a finite collection of periodic saddles $\{P_i\}_{i=1, 2, \ldots, m}$ which is $\varepsilon/2$ -dense in $H(P_f, f)$, and their Lyapunov exponents $\chi^k(P_i)$ are $\varepsilon$-close to $\chi^k(P_f)$. Moreover, as $f\in \cG$, by
{\bf R2)}
these saddles
are homoclinically related. This implies that there exists a locally maximal transitive hyperbolic set $\Lambda_f$ with $P_i\in \Lambda_f$, for all $i=1,\dots,m$. Then there is a periodic saddle $Y_f\in \Lambda_f$ that spends an arbitrary large portion of time in an arbitrary small neighborhood of the initial collection of saddles, and is $\varepsilon/2$-dense in $\Lambda_f$. This implies that the orbit of $Y_f$ is $\varepsilon$-dense in $H(P_f, f)$ and its Lyapunov exponent $\chi^k(Y_f)$ is close to the Lyapunov exponent $\chi^k(P_f)$, completing the proof of the lemma.\end{proof}
\end{description}

\subsection{Heterodimensional cycles}

In this section, we state results that allow us to generate
heterodimesional cycles for saddles in a non-hyperbolic homoclinic classes. Roughly speaking, we need to be able to create a cycle associated with a given pair of homoclinically linked saddles of
different indices (Proposition \ref{p.cycles}); to produce a locally dense set of diffeomorphisms with these type of cycles (Proposition \ref{p.simpleanddense}); and, finally, to use these cycles to generate periodic saddles with some special properties (Proposition \ref{p.bdf}) that will allow to use Proposition \ref{p.firstmainproposition} later to obtain non-hyperbolic ergodic measures with large support. 

The following result, which  is a consequence of the Connecting Lemma in
\cite{H}, can be found in \cite[Proposition 3.5]{DG}.
Recall that two saddles $P$ and $Q$ of different indices have a {\emph{heterodimensional cycle}}
if its invariant manifolds meet cyclically, $W^{\st}(P)\cap W^{\ut}(Q)\ne \emptyset$ and
$W^{\ut}(P)\cap W^{\st}(Q)\ne \emptyset$. The cycle has {\emph{coindex one}} if $\st$-index$(P)=
\st$-index$Q\pm 1$.

\bpro\label{p.cycles} Let $\, \cU$ be an open set of $\,\diffM$
such that there are saddles $P_f$ and $Q_f$ (depending
continuously on $f\in \cU$) with consecutive indices  which are
persistently linked in \ $\cU$. Then there is a dense subset $\cH$
of $\,\cU$ such that every diffeomorphism $f\in \cH$ has a coindex
one heterodimensional cycle associated to saddles $A_f\in
\text{{\rm Per}}_{\RR}(H(P_f,f))$ and $B_f\in
\text{{\rm Per}}_{\RR}(H(Q_f,f))$.
\label{p.bdhayashi} \epro

The following statements are straightforward  reformulations of Propositions 4.3 and 4.5 from \cite{DG} better suited for our case.

\bpro \label{p.simpleanddense} Let
$f$ have 
 a cycle associated to saddles $A_f$ and $B_f$ such that
\begin{itemize}
\item
$A_f$ and $B_f$ have real multipliers and $\mbox{\rm
$\st$-index\,}(A_f)=\mbox{\rm $\st$-index\,}(B_f)+1$;
\item
$A_f$ is $\st$-biaccumulated and $B_f$ is $\ut$-biaccumulated;
\end{itemize}
 Then arbitrarily $C^1$-close to $f$ there
are an open set $\mathcal{E}\subset \text{\rm Diff}^{\,1}(M)$ and
dense subset $\mathcal{D}\subset \mathcal{E}$ such that every
$g\in \mathcal{D}$ has a 
cycle associated with $A_g$
and $B_g$. \epro


\bpro{\em{(\cite[Proposition 4.5]{DG}})} \label{p.bdf} Let $f$ be
a diffeomorphism with a heterodimensional 
 cycle associated to saddles
$A_f$ and $B_f$ such that
\begin{description}
\item{{\bf (i)}}
the saddles $A_f$ and $B_f$ have real multipliers;
\item{{\bf (ii)}}
$\st$-$\text{\rm index}(A_f)=i+1$ and $\st$-$\text{\rm
index}(B_f)=i$;
\item{{\bf (iii)}}
$A_f$ is $\st$-biaccumulated, and  $B_f$ is
$\ut$-biaccu\-mu\-la\-ted.
\end{description}

Fix  neighborhoods $U_B$ of the orbit of $B_f$ and   $U_A$ of the
orbit of $A_f$. Then there are  sequences of natural numbers $\ell_k,m_k$, that
tend to infinity as $k\to \infty$, and a sequence of
diffeomorphisms $f_k$,
 $f_k\to f$ as $k\to \infty$, such that $f_k$ coincides with $f$ along the orbits of $A_{f}$ and $B_{f}$,
 and has a hyperbolic saddle $R_{k}$ 
 having real multipliers  with the following properties:
\begin{description}
\item{{\bf (1)}}
 the orbit of the saddle $R_{k}$ spends a fixed
number $t_{(a,b)}$ (independent of $k$) of iterates
  to go from  $U_B$ to $U_A$, then it remains $\ell_k\,\pi(A_f)$ iterates in  $U_A$, then it takes
a fixed number of iterates $t_{(b,a)}$ (independent of $k$) to go
from $U_A$ to $U_B$, and finally it remains $m_k\,\pi(B_f)$
iterates in $U_B$. In particular, there is a  constant $t\in \NN$
 independent of $k$ such that the period of $R_k$ is
$ \pi(R_{k})=m_k\, \pi(B_{f})+ \ell_k \, \pi(A_{f})+ t; $
\item{{\bf (2)}}
there is a constant $\const>0$ independent of $k$ such that the
central multiplier of $R_{k}$ 
satisfies $ {\const}^{-1}<|\la^{i+1} (R_{k})|<\const.$
\end{description}
Suppose also that $|\la^{i+1}(A_{f})|\in \left(\frac{1}{\sqrt{2}}, 1\right)$\footnote{In \cite{DG} here just assumed that the $(i+1)$-th multiplier $\la^{i+1}(A_{f})$ of $A_f$ is close to one. To be more accurate, here we replaced it by the condition that we actually needed for the proof, see equation (8) in \cite{DG}.}. 
 Then
 \begin{description}
\item{{\bf (3)}}
 $R_{k}$ has the same index as $B_{f}$ and is
homoclinically related to $B_{f}$;
\item{{\bf (4)}}
$W^{\st}(R_k)$ intersects  $W^{\uut}(R_k)$, 
and
$W^{\uut}(R_k)$ intersects $W^{\st}(B_f)$. Moreover, these 
intersections are
 quasi-transverse;
\item{{\bf (5)}}
there is a heterodimensional 
cycle associated to $R_{k}$ and $A_{f}$.
\end{description}
\epro

\brm\label{r.delta} Notice that the condition (1) implies that
there exists a sequence $\delta_k\to 0^+$ such that
$$
\cO (A_f), \cO(B_f)\subseteq B_{\delta_k}\left(\cO(R_k)\right).
$$
Here $\cO (X)$ denotes the full orbit of the point $X$. \erm

\brm In {\em \cite{DG}} we introduced $V$-related cycles in order to guarantee that all of the sets we consider are in the domain where
a central direction field is defined. Here this property is one of the assumptions of the main result, so we do not need to focus on $V$-related cycles. Propositions \ref{p.simpleanddense}  and  \ref{p.bdf} above are therefore weaker versions of Propositions  4.3 and  4.5  from {\em \cite{DG}}.  \erm

\section{Main construction: tree of periodic points and limit measures}\label{s.final}



 In this section we use the results of previous sections to complete the proof of Theorem~\ref{t.bdg}.
Let $\mathcal{G}$ be the residual set of $\text{\rm
Diff}^{\,1}(M)$ described in Section~\ref{ss.generic}. The statement below
is a local version of Theorem~\ref{t.bdg}.

\bthm\label{t.two} Let $f\in \mathcal{G}$
have a homoclinic class $H(f)$ of $f$ such that
\begin{itemize}
\item
$H(f)$ has a dominated splitting of the form $T_{H(f)}M=E\oplus F
\oplus G$, where $F$ has dimension one,
\item
$H(f)=H(P_f, f)=H(Q_f, f),$ where the saddles $P_f$ and $Q_f$ 
have $\st$-indices $\dim(E\oplus F)=\dim (E)+1$ and \ $\dim (E)$.
\end{itemize}
Then arbitrarily $C^1$-close to $f$ there exists a $C^1$-open set
$\mathcal{Z}\subset \text{\rm Diff}^{\,1}(M)$ and a residual
subset $\mathcal{R}\subset \mathcal{Z}$ such that every $g\in
\mathcal{R}$ has a non-hyperbolic ergodic invariant measure  whose
support is the whole homoclinic class $H(P_g, g)$.
\ethm

As in \cite{DG} (see also \cite{ABCDW}), standard genericity arguments imply that
Theorem~\ref{t.bdg} follows from Theorem~\ref{t.two}.

Before we continue, let us remind that existence of the dominated splitting is preserved by small perturbations.

\blm[Dominated splittings, Appendix
B.1.1 in \cite{BDV2}]\label{l.dominated}  Assume  that the homoclinic class
$H(P_f,f)$ has a dominated splitting $T_{H(P_f,f)}M=E_f\oplus F_f
\oplus G_f$. There are  neighborhoods $V$ of
$H(P_f,f)$ and $\,\cU_f$ of $f$ such that for every
$g\in \cU_f$, the maximal invariant set $\La_g(V)$ has a dominated splitting
$T_{\La_g(V)}M=E_g\oplus F_g \oplus G_g$, where $\dim K_f=\dim
K_g$, $K=E,F,G$. With a slight abuse of notation,  we
will omit the dependence on $g$ of this splitting.
\elm

\subsection{Construction of the sequences of periodic orbits}

Consider a homoclinic class $H(f)=H(P_f, f)=H(Q_f, f)$ of $f\in \cG$
that satisfy the conditions of Theorem \ref{t.two}.

The genericity hypotheses {\bf R2)}--{\bf R3)} imply that the homoclinic class $H(f)$  
contains two saddles $A_f$ and $B_f$
 with real multipliers such that
 $$
 \st\text{\rm -index}(A_f)=\st\text{\rm -index}(P_f),\ \ \st\text{\rm -index}(B_f)=\st\text{\rm -index}(Q_f),
 $$ and
 $$H(f)=H(A_f, f)=H(B_f, f)=H(P_f,f)=H(Q_f, f).$$
 The homoclinic class $H(f)$ has a splitting $E\oplus F \oplus
G$,  where $F$ is one-dimensional. 
We can assume that the saddle $A_f$ is $\st$-biaccumulated, the
saddle $B_f$ is $\ut$-biaccumulated, and that the central
multiplier of $A_f$ is close to one (it is enough to have
$|\lambda^{F}(A_f)|\in (0.9, 1)$), see generic conditions {\bf
R3)}--{\bf R4)}.

Recall that $\chi^{F}(A_f)$ denotes the
central Lyapunov exponent of saddle
$A_f$ corresponding to the central direction
$F$. 
Fix a constant $C$ such
that
\begin{equation}
\label{e.C} C>\displaystyle{\frac{32}{\nu}}, \ \text{\rm where}\ \nu=|\chi^F(A_f,f)|.
\end{equation}

\blm\label{l.z}
Let $f$ satisfy the conditions of Theorem \ref{t.two} and constants $C, \nu>0$ be as in \eqref{e.C}. Given $\epsilon>0$, arbitrary close to $f$ there exists an open set $\cZ\subset \diffM$  with the following properties.
\begin{description}
\item{{\bf 1)}} For every $g\in \cZ \cap \cG$
$$
H(f)\subset B_{{\epsilon/2}}(H(P_g, g)), \ \ H(P_g, g)\subset B_{{\epsilon/2}}(H(f)).
$$
\item{{\bf 2)}}
For every $g\in \cZ \cap \cG$, we have
\begin{itemize}
 \item
 $H(P_g, g)=H(Q_g, g)=H(A_g, g)=H(B_g, g)$, 
 \item there is a dominated splitting (with
one-dimensional central direction) on $H(g)=H(P_g, g)$.
\end{itemize} 

\item{{\bf 3)}}
For each  $g\in \cZ$ there is a saddle $Y_g\subset H(P_g, g)$ with real multipliers such that the following holds.
\begin{itemize}
\item the orbit of $Y_g$ is ${\epsilon}$-dense in $H(g)$,
\item
$Y_g$ and $P_g$ are homoclinically related, hence $H(Y_g,g)= H(P_g, g)$,
\item
$\chi_{Y_g}^F\in [-2\,\nu, -\nu/2]$, 
\item the saddle $Y_g$ is $\st$-biaccumulated.
\end{itemize}


\item{{\bf 4)}}
 There is  a dense in $\mathcal{Z}$
countable subset $\mathcal{D}\subset \mathcal{Z}$ such that every
$\phi\in \mathcal{D}$ has a cycle associated with the saddles
$Y_{\phi}$ and $B_{\phi}$.

\end{description}
\elm
\begin{proof}
 Since $f\in \cG$, 
by property {\bf R1)}, continuous dependence of homoclinic classes
of diffeomorphisms in $\mathcal{G}$, 
there is an open set  $\mathcal{W}\subset \text{\rm Diff}^1(M)$, $f\in \mathcal{W}$, such that for every $g\in \mathcal{W}\cap \cG$ we have
$$
H(f)\subset B_{{\epsilon/2}}(H(P_g, g)), \ \ H(P_g, g)\subset B_{{\epsilon/2}}(H(f)),
$$
obtaining property {\bf 1)}. 

By Lemma~\ref{l.dominated} there is a dominated splitting (with
one-dimensional central direction) on $H(g)=H(P_g, g)$ which is a
continuation of $E\oplus F\oplus G$.  In particular, such a
splitting also has one dimensional central direction. Slightly
abusing the notation, we will continue to denote this splitting by
$E\oplus F\oplus G$.
If $\cW$ is small enough then for every $g\in \mathcal{W}\cap \cG$, we have $H(P_g, g)=H(Q_g, g)=H(A_g, g)=H(B_g, g)$, see  {\bf
R2)}. This gives {\bf 2)}.

Due to {\bf R5)} there is $g$ arbitrary $C^1$ close to $f$ and a saddle $Y_g\subset H(P_g, g)$ with real multipliers such that its orbit is ${\epsilon/3}$-dense in $H(g)$,
$H(Y_g,g)= H(P_g, g)$,
$\chi_{Y_g}^F\in [-2\,\nu, -\nu/2]$, and which is $\st$-biaccumulated, hence we have {\bf 3)}. 

Finally, Proposition~\ref{p.bdhayashi} implies that $C^ 1$-arbitrarily close to $g$ there is $h$
 with a cycle
corresponding to the  saddles $Y_h$ and $B_h$.
Therefore, by Proposition~\ref{p.simpleanddense}, there is an open set
$\mathcal{Z}\subset \mathcal{W}$ and a dense in $\mathcal{Z}$
countable subset $\mathcal{D}\subset \mathcal{Z}$ such that every
$\phi\in \mathcal{D}$ has a cycle associated with the saddles
$Y_{\phi}$ and $B_{\phi}$. This completes the proof of Lemma \ref{l.z}. 
\end{proof}



 Next proposition can be considered as an improved version of Proposition 5.3 from \cite{DG}. Notice that  properties {\bf
Z3)} and {\bf Z4)} allow to ``spread" the constructed sequence of periodic orbits
around the whole homoclinic class.

\bpro\label{p.conditions} Let $f$ satisfy the conditions of Theorem \ref{t.two} and constants $C, \nu>0$ be as in \eqref{e.C}.
Fix a decreasing sequence
$\{\epsilon_N\}_{N\in \NN}$ of  positive numbers,
$\epsilon_N\to 0^ +$ as $N\to \infty$.
Arbitrarily $C^1$-close to $f$ there exists 
a nested sequence of  open sets
$$
\ldots \mathcal{Z}_N\subset \mathcal{Z}_{N-1}\subset \ldots \mathcal{Z}_2\subset \mathcal{Z}_1\subset \text{\rm Diff}^{\,1}(M)
$$
such that the following holds.
\begin{description}
\item{{\bf Z1)}}
A set $\mathcal{Z}_N$ is a dense open subset of $\mathcal{Z}_{N-1}$.
\item{{\bf Z2)}}
Every diffeomorphism $g\in \mathcal{Z}_{N}$ has a
finite sequence of periodic saddles homoclinically related to
$B_g$ (thus of the same index as $B_g$)
$$
\{Q_{1, g}, Q_{2, g},\ldots , Q_{N, g}\}\subset
H(B_g, g)
$$
having
 real multipliers, satisfying the  $\ut$-biaccumulation property,
and of growing periods, $\pi(B_g)=\pi (Q_{1, g})<\ldots
<\pi(Q_{N, g})$. Moreover, saddles $\{Q_{1, g},
Q_{2, g},\ldots , Q_{N, g}\}$ depend
continuously on $g$ when $g$ varies over $\mathcal{Z}_{N}$.
\item{{\bf Z3)}}
For  any  ${N}\in \mathbb{N}$ there exists a countable
dense subset  $\mathcal{D}_{{N}}\subset
\mathcal{Z}_{{N}}$ such that  every $g\in
\mathcal{D}_{{N}}$ has a saddle $Y_{{N},g}$ depending continuously on $g$ such that
\begin{itemize}
\item
 the $g$-orbit
of $Y_{N,g}$ is $\epsilon_{N+1}$-dense in
$H(B_g, g)$,
\item
$\chi^F(Y_{{N},g})\in [-2\,\nu,- \nu/2]$,
\item
there is a heterodimensional cycle associated to the
saddles $Y_{{N},g}$ and $Q_{N, g}$.
\end{itemize}
\item{{\bf Z4)}}
There is a sequence of locally constant functions $\gamma_{{N}}: \mathcal{Z}_{{N}}\to (0, +\infty)$ such that for any $N\in \mathbb{N}$  and any $g\in
\mathcal{Z}_{{N+1}}$ the orbit of $Q_{{N+1}, g}$ is a $(\gamma_{{N}}(g), 1- C\,
|\chi^{F}(Q_{N, g})|)$-good approximation  of the
orbit of $Q_{N, g}$ (recall
Definition~\ref{d.goodaprox}), where $C$ is the constant in
\eqref{e.C}. Moreover,
$$
{{B}}_{\epsilon_{N+1}} (\cO (Q_{{N+1}, g})) \supset
\cO ( Y_{N, g}).
$$
\item{{\bf Z5)}}
Take any $g\in \mathcal{Z}_{{N}}$. Let $d_k, 1\le
k\le N$,  be the minimal distance between the points of the
$g$-orbit of $Q_{{k}, g}$. Then
$$
\gamma_{{N}}(g)< \frac{\min_{1\le k\le N} d_k}{3\cdot 2^N}. 
$$
\item{{\bf Z6)}}
For each $N\in \mathbb{N}$ and every $g\in \mathcal{Z}_{{N+1}}$
$$
|\chi^{F}(Q_{{N+1}, g})|< \frac{1}{2}\,
|\chi^{F}(Q_{{N}, g})|.
$$
\end{description}
\epro

\begin{proof}
We prove Proposition~\ref{p.conditions} by induction.

As a set $\mathcal{Z}_1$ we take the set $\mathcal{Z}$ from Lemma \ref{l.z} (where we take $\epsilon=\epsilon_2$). Taking $Q_{1,g}=B_g$, $\mathcal{D}_1=\mathcal{D}$, and $Y_{1,g}=Y_g$, we see that properties {\bf Z2)} and {\bf Z3)} for $\mathcal{Z}_1$ are satisfied. This form the base of induction. Notice that here we do not need to check the properties {\bf Z4)} - {\bf Z6)}.

Now assume that the sets $\mathcal{Z}_1 \supset \mathcal{Z}_2 \supset \ldots \mathcal{Z}_N$ together with periodic orbits $\{Q_{1, g}, \ldots Q_{N, g}\}$, $\{Y_{1, g}, \ldots Y_{N, g}\}$, functions $\gamma_i:\mathcal{Z}_i\to (0, +\infty)$, and sets $\mathcal{D}_i$  are constructed,  $ i=1, \ldots, N$. Let us construct the set $\mathcal{Z}_{N+1}\subset \mathcal{Z}_N$.

 Notice that it is enough to construct an open dense subset of each connected component of $\mathcal{Z}_N$ (together with periodic orbits $Q_{N+1, g}, Y_{N+1, g}$, function $\gamma_{N+1}$, and a dense subset $\mathcal{D}_{N+1}$) that satisfies conditions ${\bf Z2)}$ - ${\bf Z6)}$.  Take one of the connected components of the set $\mathcal{Z}_N$, denote it by $\mathcal{Z}^*$.

Let us
enumerate diffeomorphisms from  $\mathcal{D}_N\cap \mathcal{Z}^*=\{g_i\}_{i\in
\mathbb{N}}$. Take one of these diffeomorphisms, say $g_i\in
\mathcal{D}_N\cap \mathcal{Z}^*$.

Proposition~\ref{p.bdf} (applied to the saddles $Q_{N,g_i}$ and $Y_{N,g_i}$)
 allows to obtain a sequence of
diffeomorphisms $g_{ik}$, $g_{ik}\to g_i$ as $k\to \infty$, such
that each diffeomorphism $g_{ik}$ has a periodic saddle $S_{ik}$
with real multipliers (denoted by $R_k$ in
Proposition~\ref{p.bdf}), having the following
properties:
\begin{description}
\item{\bf S1)}
the saddle $S_{ik}$ is homoclinically related to $Q_{N, g_{ik}}$,
thus has the same index as $B_{g_i}$;
\item{\bf S2)}
for some constant $\const$ that does not depend on $k$ and for
each $k\in \mathbb{N}$ we have $1<|\lambda^F(S_{ik})|<\const$;
\item{\bf S3)}
the map $g_{ik}$ has a heterodimensional cycle associated to
$S_{ik}$ and $Y_{N, g_{ik}}$;
\item{\bf S4)}
take a real number $\gamma_N=\gamma_N(\mathcal{Z}^*)>0$ small enough (the choice of $\ga_N$ will be clear from the construction below)
and
$\ga_N$-neighborhoods $U_{\cO(Y_{N, g_i})}$ and
 $U_{\cO(Q_{N, g_i})}$ of the orbits of $Y_{N, g_i}$ and $Q_{N, g_i}$, and sequences of natural numbers $\ell_k,m_k$ that tend to
infinity as $k\to \infty$, such that under the iterates of
$g_{ik}$ it takes to the saddle $S_{ik}$  a fixed number of iterates 
 (independent of $k$) to go from  $U_{\cO(Y_{N, g_i})}$
 to $U_{\cO(Q_{N, g_i})}$, then it remains $\ell_k\,\pi(Y_{N, g_i})$ iterates in  $U_{\cO(Y_{N, g_i})}$,
then it needs a fixed number of iterates 
to go from $U_{\cO(Y_{N, g_i})}$ to $U_{\cO(Q_{N, g_i})}$, and finally it remains
$m_k\,\pi(Q_{N, g_i})$ iterates in $U_{\cO(Q_{N, g_i})}$.
In particular,
there is a  constant $t\in \NN$
 independent of $k$ such that
$$
\pi(S_{ik})=m_k\, \pi(Q_{N, g_i})+ \ell_k \, \pi(Y_{N, g_i})+ t.
$$
\end{description}
Exactly as in \cite{DG}, properties (3) and (4) from Proposition~\ref{p.bdf}
guarantee that making an arbitrary small perturbation of $g_{ik}$
(preserving properties {\bf S1)} - {\bf S4)}) we also have:
\begin{description}
\item{\bf S5)}
the saddle $S_{ik}$ has the $\ut$-biaccumulation property.
\end{description}

We need the following  quantitative estimates:
\blm
[Lemma 5.4 from \cite{DG}]
\label{l.inequalities} For every large $k\in \mathbb{N}$ the
saddles $S_{ik}$ satisfy:
\begin{equation}\label{e.1}
   0<\chi^F(S_{ik})<\frac{1}{2}\chi^F(Q_{N, g_{ik}}),
 \end{equation}
 \begin{equation}\label{e.2}
\frac{m_k \, \pi (Q_{N, g_i})}{\pi(S_{i,k})}=
    \frac{m_k\pi(Q_{N, g_i})}{m_k\pi(Q_{N, g_i})+\ell_k\pi(Y_{N, g_i})+t}>1-C\chi^F(Q_{N, g_{ik}}).
\end{equation}
\elm

Consider now the set
$$
\mathcal{D'}=\left\{ g_{ik}\  |\  S_{ik} \text{\rm \ \ satisfies
conditions (\ref{e.1}), (\ref{e.2}), $i\ge 0$ and $k\ge 0$} \right\} \subset \mathcal{Z}^*.
$$
By construction, the set $\mathcal{D'}$ is a countable dense
subset of $\mathcal{Z}^*$. Let us enumerate the elements of
$\mathcal{D'}=\{h_{n}\}_{n\in \NN}$.  Let us also denote by
$Q_{N+1, g}(n)$ the continuation of the periodic saddle $S_{ik}$ of the map $h_{n}\equiv
g_{ik}$. In particular, $Q_{N+1, g_{ik}}(n)=S_{ik}$.

Now an application of Lemma \ref{l.perdense} and  Proposition~\ref{p.simpleanddense} implies that  for each $h_{n}$  there is an open set $\mathcal{U}_n$, which is
$\frac{1}{n}$-close to
 $h_{n}$, and for each $g\in \mathcal{U}_n\cap \mathcal{\cG}$ there is  a saddle $Y_{N+1,g}(n)$ such that
   \begin{description}
   \item{{\bf Y1)}}
    the saddle  $Y_{N+1,g}(n)$  is in $H(g)=H(P_g,g)$ and its orbit is $\epsilon_{N+2}$-dense in $H(g)$,
   \item{{\bf Y2)}}
   the saddles
   $Q_{N+1,g}(n)$ and $Y_{N+1,g}(n)$ have different indices,
   \item{{\bf Y3)}}
   the Lyapunov exponent $\chi^F(Y_{N+1,g}(n))$ is close to
   $(-\nu)$, and
   \item{{\bf Y4)}}
 there is a dense subset $\widetilde{\mathcal{D}}{(n)}$ of maps in $\mathcal{U}_n$ having a heterodimensional cycle between $Q_{N+1,g}(n)$ and $Y_{N+1,g}(n)$.
   \end{description}

%

        We can take
    $\mathcal{U}_{n}$   small enough and close enough to $h_n$ to guarantee that for every $g\in \mathcal{U}_{n}$ one has
    $0<\chi^F(Q_{N+1, g}(n))<\frac{1}{2}\chi^F(Q_{N,g})$. Indeed, due to
    (\ref{e.1}) this inequality holds for $h_{n}$. Since
    Lyapunov exponents of a hyperbolic saddle depend continuously on
    a diffeomorphism, the inequality holds also for all $g$ sufficiently $C^1$-close to
    $h_{n}$.

Let us now define inductively
$$
\mathcal{Z}(1)=\mathcal{U}_{1}, \
\mathcal{Z}(2)=\mathcal{U}_{2}\backslash \overline{\mathcal{Z}(1)},
\ldots, \mathcal{Z}({n})=\mathcal{U}_{n}\backslash
\overline{\mathcal{Z}({n-1})}, \ldots
$$
and
$$
\mathcal{D}({n})=\widetilde{\mathcal{D}}{(n)}\cap
\mathcal{Z}({n}), n\in \NN.
$$
Also we set
$$
\mathcal{Z}_{N+1}^*=\bigcup_{{n\in \NN}}\mathcal{Z}({n}), \qquad
 \mathcal{D}_{N+1}^*=\bigcup_{n\in \NN}\widetilde{\mathcal{D}}{(n)}.
$$
Besides, for each $g\in \mathcal{Z}({n})$ we set
$Q_{N+1,g}=Q_{N+1,g}(n)$ and $Y_{N+1,g}=Y_{N+1,g}(n)$. Finally, we define $\mathcal{Z}_{N+1}$ and $\mathcal{D}_{N+1}$ as the union of constructed sets $\mathcal{Z}_{N+1}^*$ and $\mathcal{D}_{N+1}^*$ over all connected components of $\mathcal{Z}_{N}$.
We claim that the  set
$\mathcal{Z}_{N+1}$ satisfies the required
properties {\bf Z1)} - {\bf Z6)}.
\begin{itemize}
\item
Since the set $\{h_{n}\}_{n\in \NN}$ is dense in 
$\mathcal{Z}^*$ (the chosen connected component of $\mathcal{Z}_N$), the union $\cup_{n\in \NN}\mathcal{U}_{n}$ is
dense in $\mathcal{Z}^*$, and hence $\cup_{n\in
\NN}\mathcal{Z}({n})$ is also dense in $\mathcal{Z}^*$, so {\bf
Z1)} holds.
\item
For each $g\in \mathcal{Z}({n})$ the saddle $Q_{N+1, g}(n)$ is homoclinically related to $B_g$ (recall {\bf S1)}), has
real multipliers, $\ut$-biaccumulation property, and
$\pi(Q_{N, g})<\pi(Q_{N+1, g}(n))$, so {\bf Z2)} holds.
\item
The sets $\widetilde{\mathcal{D}}({n})\subset \mathcal{U}_{n}$ and the saddles $Y_{N+1, g}({n})$ were
constructed to satisfy {\bf Z3)}. More precisely, the estimate of the Lyapunov exponent
$\chi^F(Y_{N+1, g}({n}))$ follows from the choice of $Y_{N+1, g}({n})$ in {\bf Y3)}, the $\epsilon_{N+2}$-density of the orbit of
$Y_{N+1, g}({n})$ in $H(g)$ follows from {\bf Y1)}, and the existence of the heterodimensional cycle follows from {\bf Y4)}.
\item
Take $g\in \mathcal{Z}_{N+1}$, and denote by $\Gamma$ the part of
the orbit of $Q_{N+1, g}$ that belongs to the neighborhood $U_{Q_{N,g}}$.
Define the projection
$$
\rho:\Gamma\to \mathcal{O}(Q_{N,g}), \quad \rho(x)=\{\text{\rm the
point of $\mathcal{O}(Q_{N,g})$ nearest to $x$}\}.
$$
By construction (recall {\bf S4)}),
$$
\#\Gamma=m_k\pi(Q_{N,g})\quad \mbox{and} \quad
\#(\mathcal{O}(Q_{N+1, g}))=m_k\pi(Q_{N,g})+\ell_k\pi(Y_{N,g})+t.
$$
Recall that
 here $k$
and $n$ are related due to the enumeration $h_{n}=g_{ik}$;
notice that in fact integers $m_k$, $\ell_k$, and $t$ depend also
on the index $i$, but our notations do not reflect this
dependence. Now the first part of {\bf Z4)} follows from the inequality (\ref{e.2}).
The second part,
${{B}}_{\epsilon_{N+2}} (\cO (Q_{N+1, g})) \supset
\cO ( Y_{N,g})$, follows from the fact that $\gamma_N$ can be chosen arbitrarily small,
in particular, smaller than $\epsilon_{N+2}$.
\item
Condition {\bf Z5)} is obtained by  choosing
sufficiently small $\gamma_N>0$.
\item
The last property {\bf Z6)} follows directly from inequality
(\ref{e.1}).
\end{itemize}

We completed the step of induction, and, thus, proved  Proposition~\ref{p.conditions}.
\end{proof}

\subsection{Infinite sequences of periodic orbits}

Here we conclude the proof of Theorem \ref{t.two} and, therefore, of the main result.

\bpro\label{p.localperiod}
 Let $f$ satisfy the conditions of Theorem \ref{t.two}. Arbitrary close to $f$ there exists an open set $\cZ\subset \diffM$ such that generic diffeomorphism $g$ from $\mathcal{Z}$ have a
sequence of periodic saddles 
which satisfies the
assumptions of Proposition~\ref{p.firstmainproposition}, belongs
to the  homoclinic class $H(P_g, g)$, and the union
of their orbits is dense in the homoclinic class $H(P_g,g)$.
 \epro

Note that Propositions~\ref{p.localperiod} and
\ref{p.firstmainproposition} give non-hyperbolic
ergodic measures supported on the whole homoclinic class for generic diffeomorphisms from $\mathcal{Z}$,
and, thus, Theorem~\ref{t.two}.

\medskip

\begin{proof}
Apply Proposition \ref{p.conditions}, and set $\mathcal{Z}=\mathcal{Z}_1$. Due to Property {\bf Z1)}, for any $N\in \mathbb{N}$ the set
${\mathcal{Z}}_{N}$ is an open and
dense subset of $\mathcal{Z}$. Therefore the intersection
$$
\mathcal{R}=\cG\cap\left(\bigcap_{N\in \mathbb{N}}{\mathcal{Z}}_N\right).
$$
 is a residual subset of $\mathcal{Z}$.

 Take any $g\in \mathcal{R}$.  Due to {\bf Z2)}, for the
diffeomorphism $g\in \mathcal{R}$ a sequence of periodic points
$$\{B_g, Q_{1, g}, Q_{2, g}\ldots ,
Q_{N, g}, \ldots \}\subset H(P_g, g)$$ is
well defined. We claim that this sequence satisfies the
assumptions of Proposition~\ref{p.firstmainproposition}.

Indeed, assumption {\bf 1)} holds since $g\in \cG$ and due to the choice of $\cZ$ 
(existence of a one dimensional center
direction for $\Delta=H(P_g,g)$), and {\bf 2)} follows from 
property {\bf Z2)}. Assumptions {\bf
3)}, {\bf 4)}, and {\bf 5)} follow from {\bf Z4)},   {\bf Z5)},
and {\bf Z6)}, respectively.

Finally, since $\epsilon_N\to 0$ as $N\to \infty$, from {\bf Z3)}, {\bf Z4)} we have
that the nonhyperbolic
ergodic measure given by Proposition~\ref{p.firstmainproposition} satisfy
\begin{equation}
\label{e.formula}
\mathrm{supp} (\mu)=
 \bigcap_{N\in \NN} \quad \left( \overline{
\bigcup_{ K\ge N}
 \cO
(Q_{K, g})         } \right) =H(P_g, g).
\end{equation}
This proves
Proposition~\ref{p.localperiod}.
\end{proof}

\vskip 1cm

\flushleft{\bf Christian Bonatti}
 \ \ (bonatti@u-bourgogne.fr)\\
Institut de Math\'ematiques de Bourgogne\\ B.P. 47 870\\
21078 Dijon Cedex \\ France
\medskip

\medskip

\flushleft{\bf Lorenzo J. D\'\i az}
 \ \ (lodiaz@mat.puc-rio.br)\\
Departamento de  Matem\'{a}tica, PUC-Rio \\ Marqu\^{e}s de S. Vicente 225\\
22453-900 Rio de Janeiro RJ \\ Brazil

\medskip

\flushleft
{\bf Anton Gorodetski}  \ \  (asgor@math.uci.edu)\\
 Department of Mathematics\\
University of California, Irvine\\
 Irvine, CA 92697\\
USA 

\end{document}